\documentclass[12pt]{article}
\usepackage{amsmath,amsfonts,amsthm,amscd}
\newtheorem{theorem}{Theorem}

\newtheorem{proposition}{Proposition}

\newtheorem{corollary}{Corollary}

\numberwithin{equation}{section}
\numberwithin{theorem}{section}
\numberwithin{proposition}{section}
\numberwithin{lemma}{section}
\numberwithin{claim}{section}
\numberwithin{corollary}{section}
\newcommand{\cpn}{\ensuremath{\mathbf{C}P^{N}}}

\newcommand{\dl}{\ensuremath{\partial}}
\newcommand{\dlb}{\ensuremath{\overline{\partial}}}

\newcommand{\bull}{\ensuremath{{}\bullet{}}}

\newcommand{\gc}{\ensuremath{G^{\mathbb{C}}}}

\input{diagrams.tex}

\newarrow{to}---->
\newarrow{mto}|--->
\newarrowtail{<=}\Rightarrow\Leftarrow\Uparrow\Downarrow
\newarrow{bboth}{<=}==={=>}
\newarrow{inject}{hooka}---{vee}
\newarrow{surject}----{>>}
\begin{document}
\title{Algebraic and Analytic K-Stability}

\author{{Sean T. Paul}\thanks{Research supported in
part by an NSF Postdoctoral Fellowship}
\and{Gang Tian}\thanks {Massachusetts Institute of Technology}}

\date{May 26, 2004}
\maketitle
\vspace{-5mm}
\begin{abstract}
In this note we identify the leading terms of the (reduced) K-energy map with a universal linear combination of the principal
and subdominant coefficients of the weight of the $mth$ Hilbert point. This shows that the weight $F_{1}(\lambda;X)$ introduced
by Donaldson in [SKD02] is just the weight of the CM-polarisation.The equivalence between the CM-(semi)stability and the K-(semi)
stability follows from this. 
Also, using our previous work, we are able to describe
this subdominant coefficient in terms of the weights of some generalised Chow forms, under a multiplicity free hypothesis on the degeneration. This is accomplished by introducing a parameter dependent
lift of the CM-polarisation, and letting this parameter tend to infinity. This could be thought of as a ``quantized'' version of the virtual bundle introduced in [Tian94].
\end{abstract}
\newpage
\section{Introduction and Motivation}
One of the central problems in complex differential geometry is to
find necessary and sufficient conditions for the existence of
``canonical'' metrics within a given K\"{a}hler class. In the early
80's, E. Calabi introduced the notion of extremal metrics. Most
extremal metrics are in fact K\"ahler metrics of constant scalar
curvature, however it is still very difficult to find a  K\"ahler
metric of constant scalar curvature in a general K\"ahler class.
In the late 80's, Yau conjectured that the existence of K\"ahler-Einstein metrics with positive scalar
curvature should be related to the stability (in the sense of Mumfords' G.I.T.) of the underlying algebraic manifold.
In [Tian97], the second author introduced the notion of \emph{CM stability} and
\emph{K-stability} and proved that the existence of K\"ahler-Einstein metrics implies those stabilities.
The arguments in [Tian97] also provided strong evidence that these two stabilities are closely related and should be
equivalent to existence of K\"ahler metrics of constant scalar curvature (also see [Tian00]).
In [Don01], S.K. Donaldson proved that the existence of K\"ahler metrics of constant scalar curvature implies
the Chow or Hilbert (semi)stability. In the last few years, there have been
many exciting works on the geometric stability of projective manifolds
([Don02], [Ross-Thomas],[PhSt02] [Futaki04], [Mabuchi04]). However, it has not been clear how the
CM- or K-stabilities are related to the standard notions of
stability coming from G.I.T. e.g., Chow and Hilbert stability. The
purpose of this paper is to examine precise connections among some of these stabilities.
In order to better understand our results, it is instructive to
first recall the well-known picture for holomorphic vector
bundles. It can be regarded as a ``linearised'' version of the
problem we are studying.

Let $\mathcal{E}$ be a coherent sheaf on a polarized manifold $(X,\mathcal{L})$.
As usual we define the Hilbert polynomial of $\mathcal{E}$ relative to $\mathcal{L}$
to be
\begin{align*}
P^{\mathcal{E}}(m):= h^{0}(\mathcal{E}(m))
\end{align*}

There are several notions of stability for vector bundles. These are not at all equivalent, for our purposes
there are two kinds of stability which will be singled out for attention.
Let $\mathcal{E}$ be a coherent sheaf on $(X,\mathcal{L})$ with Hilbert polynomial
$P^{\mathcal{E}}(m)$, then $\mathcal{E}$ is (Gieseker)\emph{stable} iff for every coherent subsheaf $\mathcal{F}$ of $\mathcal{E}$
we have the inequality
\begin{align*}
\frac{P^{\mathcal{F}}(m)}{\mbox{rnk}{\mathcal{F}}} < \frac{P^{\mathcal{E}}(m)}{\mbox{rnk}{\mathcal{E}}}.
\end{align*}
The semistability allows possible equality.
On the other hand,  we say that $\mathcal{E}$ is (Mumford) \emph{stable} iff for every coherent subsheaf $\mathcal{F}$ of $\mathcal{E}$
we have the inequality
\begin{align*}
\frac{\mbox{deg}_{\omega}(\mathcal{F})}{\mbox{rnk}\mathcal{F}}<
\frac{\mbox{deg}_{\omega}(\mathcal{E})}{\mbox{rnk}\mathcal{E}}.
\end{align*}
Again the semi-stability allows possible equality. It is not hard to see that Mumford stability implies Gieseker stability.
If ${\mathcal E}$ is an irreducible vector bundle over a polarized K\"ahler manifold
$(X,\mathcal{L})$, then the Donaldson-Uhlenbeck-Yau theorem states that
$\mathcal{E}$ admits a Hermitian-Einstein metric with respect to $\omega=c_1({\mathcal L})$ iff
$\mathcal{E}$ is Mumford stable with respect to $\omega$. This was first conjectured by Hitchin-Kobayashi. Here by a \emph{Hermitian-Einstein} metric, we mean a Hermitian metric $h$ on ${\mathcal E}$ whose curvature tensor
satisfies
\begin{align*}
{\rm Tr}_{\omega}F_{h}=\lambda \mbox{Id}_{\mathcal{E}} ,\qquad \lambda \ \mbox{a constant}
\end{align*}
This closely parallels the relationship between Stability of the Hilbert point
and the K-Stability, which are the relevant concepts
for the constant scalar curvature problem.

Now let us describe our main results. First we recall the notion
of CM stability introduced in [Tian94] (also see [Tian97],
[Tian00]). Let $\pi: \mathfrak{X} \mapsto \mathcal{H}$ be a $\gc$
equivariant morphism between projective (connected) schemes over
$\mathbb{C}$ satisfying :

\begin{tabular}{ll}\\
1) The scheme $\mathfrak{X} \subset \mathcal{H}\times \cpn$ is a
family of subschemes of dimension $n$, where\\ the action of $\gc$
on $\mathfrak{X}$ is induced by the standard action on $\cpn$ and
$\pi=\pi_1|_\mathfrak{X}$.
\\
\\
2) The family $\mathfrak{X}$ is \emph{flat}, i.e., there is a
numerical polynomial $P$ such that for\\ every $z\in \mathcal{H}$
and all $m>>0$, we have
\end{tabular}
\begin{align}
h^{0}(\pi_{1}^{-1}(z),\pi_{2}^{*}{\mathcal
O}_{\mathbb{P}^{N}}(m))= P(m),
\end{align}

Consider the virtual bundle on $\mathfrak{X}$
\begin{align}
\begin{split}
2^{n+1}\mathcal{E}:= (n+1)(\mathcal{K}^{-1}-\mathcal{K})({\mathcal L}-{\mathcal L}^{-1})^{n}-
\mu({\mathcal L}-{\mathcal L}^{-1})^{n+1},\quad \\
\end{split}
\end{align}
where $\mathcal{K}$ is the relative dualizing sheaf of the
morphism and ${\mathcal L}= \pi_{2}^{*}{\mathcal
O}_{\mathbb{P}^{N}}(1)$. Then the \emph{CM polarization} is the
$\gc$ linearized line bundle defined by
\begin{align}
\textbf{L}_{CM}^{-1}:= {\bf{det}}(R_{\pi*}^{\bull}(\mathcal{E}))
\end{align}
The idea in this paper is to \emph{change the lift} \footnote{a virtual bundle $\mathcal{E}$
such that ${\bf{det}}(R_{\pi_{1}*}^{\bull}(\mathcal{E}))^{p}=
\textbf{L}_{CM}^{-1}$ , $p\in \mathbb{Q}_{+}$.} \emph{without
changing the CM polarization}. That is, we consider the
(apparently) \emph{parameter dependent} virtual bundle
$\mathcal{E}(m)$, $m>>0$
\begin{align}
\begin{split}
\mathcal{E}(m):= & (-1)^{n}2(n+1)\{
{\mathcal L}^{m}({\mathcal O}-{\mathcal L})^{n}+mL^{m}({\mathcal O}-{\mathcal L})^{n+1}\}\\
&-(-1)^{n+1}(\mu+n(n+1)){\mathcal L}^{m}({\mathcal O}-{\mathcal
L})^{n+1}
\end{split}
\end{align}
and then similarly we introduce the following polarization on
$\mathcal{H}$
\begin{align}
\textbf{L}(m):= {\bf{det}}(R_{\pi*}^{\bull}(\mathcal{E}(m))
\end{align}

Now we come to the question of creating a family $\mathfrak{X}\rightarrow\mathcal{H}$ in which our $X\subset \cpn$ will move. This can be done, by appealing to Grothendeicks' construction of the Hilbert scheme, but this is actually more than we need. 
 
We proceed as follows. Given $X\subset\cpn$ we set
\begin{align*}
G^{\mathbb{C}}X:=\{(\sigma,y)\in G^{\mathbb{C}}\times \cpn : y\in \sigma X\}
\end{align*}
Let $\overline{G^{\mathbb{C}}}$ be the DeConcini-Procesi ``wonderful'' compactification \footnote{We would like to thank Micheal Thaddeus for bringing this to our attention}of $G^{\mathbb{C}}$. This is our $\mathcal{H}$.
It is well known that this compactification is smooth and has the crucial property that the first chern class map is injective. We will call the associated $\mathfrak{X}$ the DeConcini-Procesi family associated to $X\subset \cpn$.

Let $\mathfrak{X}:= \overline{G^{\mathbb{C}}X}$ be the closure of 
$G^{\mathbb{C}}X$ inside of $\overline{G^{\mathbb{C}}}\times \cpn$.
Then $\mathfrak{X}$ has a divisor singularity. The base point corresponding to our $X$ is just the identity element. 
\begin{theorem}
\label{cmpol} 
Let $\pi: \mathfrak{X} \mapsto \mathcal{H}$ be a $\gc$
equivariant morphism between projective (connected) schemes over
$\mathbb{C}$ satisfying the conditions 1) and 2).
Assume that the connecting homomorphism
\begin{align*}
c_{1}:\mbox{H}^{1}(\mathcal{H},{\mathcal O}^{*})\rightarrow
\mbox{H}^{2}(\mathcal{H},\mathbb{Z})
\end{align*}
is injective. 
Then there is a
$\gc$-isomorphism\footnote{Actually $\gc$ can be any reductive
semisimple affine algebraic group over $\mathbb{C}$} of line
bundles on $\mathcal{H}$
\begin{align}
\textbf{L}_{CM}^{-1} = \textbf{L}(m)
\end{align}
In particular, $\textbf{L}(m)$ is independent of $m$.
\end{theorem}

This result has many consequences. Let $(X,L)$ be a polarized
manifold such that $X\subset \cpn$ and $L={\mathcal O}_{\cpn}
(1)$. Let $\lambda$ be any 1psg. of $\gc$. Then we denote by
$a_{n+1}(\lambda)$ and $a_{n}(\lambda)$ the corresponding
coefficients of the weight of the Hilbert point of $X\subset \cpn$
relative to the 1psg. $\lambda:\mathbb{C}^{*}\rightarrow \gc$. We
define $F_{1}$\footnote{This definition is due to S. Donaldson} as
follows
\begin{align}
F_{1}(\lambda\ , X):= \frac{n!}{2d}(2a_{n}-\mu a_{n+1}),
\end{align}
where $\mu$ is the average of the scalar curvature (essentially
the coefficient of $k^{n-1}$ in the Hilbert polynomial $P$).
\begin{theorem}
Let $(X, L)$ be as above, then
\begin{align}\label{weight}
w_{CM}(\lambda,X)= 2 d (n+1) F_{1}(\lambda,X)
\end{align}
Where the weight $w_{CM}(\lambda,X)$ has been computed with respect to DeConcini-Procesi family.
\end{theorem}
This follows from Theorem \ref{cmpol} by taking the
weight of both sides in (1.6) and letting
$m\rightarrow \infty$ on the right hand side.

In [PT04] the weight of the CM polarisation was described in terms of ``double'' Chow coordinates, we refer the reader to that paper for details. Precisely, when the limit cycle has no multiple components we have
\begin{align}
w_{CM}(\lambda,X)=w(\lambda,\ f_{\mathcal{D}})-\left(2d
+\frac{\mu(X)}{n+1}-(n+2)\right)w(\lambda,\ R_{X}),
\end{align}
Under these same hypothesis we deduce
\begin{corollary}
\begin{align}
2d(n+1)F_{1}(\lambda,X)=w(\lambda,\ f_{\mathcal{D}})-\left(2d
+\frac{\mu(X)}{n+1}-(n+2)\right)w(\lambda,\ R_{X}).
\end{align}
\end{corollary}

The next result depends on the main arguments of [Tian94] (also
see [Tian97], [Tian00]). In that paper, the leading term of the
reduced K-Energy map\footnote{The reduced K-energy was not
explicitly defined in [Tian94]} was identified with the weight of
the CM polarization.

Since $X\subset \cpn$ with $L={\mathcal O}_{\cpn}(1)$, for any
$\sigma \in \gc$, we have a K\"ahler metric
$\sigma^*\omega_{FS}|_X$, where $\omega_{FS}$ is the Fubini-Study
metric on $\cpn$. All such metrics can be parametrized by $\gc$
modulo $U(N+1)$. Let $\omega$ be a fixed K\"ahler metric on $X$
with K\"ahler class $c_1(L)$, then we can write
$$\sigma^* \omega_{FS}|_X = \omega + \partial \overline{\partial}
\varphi_\sigma.$$
\begin{theorem} (\emph{Asymptotics of the reduced K-energy map})
Let $(X,L)$ be as above. There is a function

\begin{align*}
\Psi_{X}:\gc \rightarrow \mathbb{R} 
\end{align*}
depending only on the
embedding of $X$ where $-\infty \leq \Psi_{X} \leq C$ such that
\begin{align}
 d\nu_{\omega}(\varphi_{\lambda(t)})-\Psi_{X}({\lambda(t)})= 4 d
 F_{1}(\lambda,X)\log(t)+O(1),
\end{align}
where $\nu_\omega$ denotes the K-energy of Mabuchi.
\end{theorem}
The function  $\Psi_{X}$
degenerates to $-\infty$ if $X^{\lambda(0)}=\lim_{t\mapsto 0}
\lambda(t) (X)$ is non-reduced. The \emph{Reduced} K-energy is
defined to be the quantity on the left hand side of the above
equation. The crucial point is that $\Psi_X$ is \emph{bounded from
above}. This term appears when one compares the CM-stability with
the extremal behavior of the K-energy map.

The following corollary was suggested to us by Julius Ross.
\begin{corollary}
Suppose that $X$  has a K\"ahler-Einstein metric and
$F_{1}(\lambda,X)=0$, then the limit cycle $X^{\lambda(0)}$ is
multiplicity free.
\end{corollary}

It was observed in [Don02] that $F_1(\lambda, X)$ coincides with
the Calabi-Futaki invariant
$Re(\mathfrak{F}_{X^{\lambda(0)}}(\lambda'(1)))$ when the limit
cycle is \emph{smooth}. Combining the above theorems with the main
result in [DT02]\footnote{The result in [DT92] was stated for $X$
polarized by its anti-canonical bundle, however, with
modifications on notations and the same arguments, one can easily
extend the main result to general polarized manifolds.} one can
show that the same holds even for the generalized Calabi-Futaki
invariant in [DT92] when the limit cycle is normal.
\begin{corollary} If $X\subset \cpn$ and $\lambda$ is an 1-psg. in
$\gc$ with normal limit cycle $X^{\lambda(0)}$, then we have
\begin{align}
Re(\mathfrak{F}_{X^{\lambda(0)}}(\lambda'(1)))= -4F_{1}(\lambda;X)
\end{align}
\end{corollary}
In particular, using the properness of the K-energy (established
in [Tian97]), one deduces that $F_{1}(\lambda;X) < 0$ for any
1-psg. $\lambda$ when $X$ admits a K\"ahler-Einstein metric and
has no non-vanishing holomorphic vector fields.

\section{Proof of Theorem 1.1}
We begin with some preliminaries on $\gc$-linearized line bundles.
This material is standard, for further information see [Thad96], [Dolg03]. Denote by
$\mbox{Pic}^{\gc}(\mathcal{H})$ the abelian group of
$\gc$-linearized line bundles on $\mathcal{H}$. We need to study
the kernel of the forgetful homomorphism
\begin{align}
\alpha:\mbox{Pic}^{\gc}(\mathcal{H})\rightarrow \mbox{Pic}(\mathcal{H})
\end{align}
It is well known that $\mbox{Ker}(\alpha)$ parametrizes all the
$\gc$-linearizations on the \emph{trivial} line bundle ${\mathcal
O}$. Any $\gc$-linearization on ${\mathcal O}$ corresponds to 
an algebraic 1-cocycle
\begin{align*}
\Psi:\gc\times \mathcal{H}\rightarrow
\mathbb{C}^{*}:~~\Psi(gg^{'},z)=\Psi(g,g^{'}z)\Psi(g^{'},z).
\end{align*}
Since $\mathcal{H}$ is a proper variety, any non-vanishing holomorphic
function must be constant. It follows that
$\Psi(g,z)=\Psi(g)$, so
$\gc$ linearizations on $\mathcal{O}$ are given by the
characters of $\gc$. However, $\gc$ is \emph{semisimple} and so
the character group is trivial. Therefore
\begin{proposition}
A line bundle on $\mathcal{H}$ admits at most one
$\gc$-linearization.
\end{proposition}
Since line bundles on $\mathcal{H}$ are determined by their first
Chern class, we only need to show that
\begin{align*}
c_{1}(\textbf{L}(m))= c_{1}(\textbf{L}_{CM}^{-1})
\end{align*}
in order to guarantee that they are isomorphic as $\gc$ bundles.
This in turn is a very straightforward consequence of the
Grothendieck-Riemann-Roch theorem. We will need the following
identity (which the reader may easily establish)
\begin{align}
\sum_{i=0}^{n}(-1)^{i}\binom{n}{i}(m+i)^{k}=\begin{cases}&
(-1)^{n}\left(m+\frac{n}{2}\right)(n+1)!\qquad k=n+1\\& (-1)^{n}n!
\qquad k= n\\& 0 \qquad k\leq n-1
\end{cases}
\end{align}
As before we let $\mathcal{E}$ denote the virtual bundle over $\mathfrak{X}$
\begin{align}
2^{n+1}\mathcal{E}:= (n+1)(\mathcal{K}^{-1}-\mathcal{K})({\mathcal
L}-{\mathcal L}^{-1})^{n}-\mu({\mathcal L}-{\mathcal
L}^{-1})^{n+1}
\end{align}
which defines the CM polarization.

We also need the following virtual bundle $\mathcal{E}(m)$
\begin{align}
\begin{split}
\mathcal{E}(m):= & (-1)^{n}2(n+1)\{{\mathcal L}^{m}({\mathcal O}-{\mathcal L})^{n}
+m{\mathcal L}^{m}({\mathcal O}-{\mathcal L})^{n+1}\}\\
&-(-1)^{n+1}(\mu+n(n+1)){\mathcal L}^{m}( {\mathcal
O}(m)-{\mathcal L})^{n+1}
\end{split}
\end{align}

Recall that $\textbf{L}(m)$ is the determinant of the direct image
of $\mathcal{E}(m)$
\begin{align}
\textbf{L}(m):= {\bf{det}}(R_{\pi*}^{\bull}(\mathcal{E}(m))\qquad
m>>0
\end{align}

We only need calculate the Chern charcater of the corresponding lift. It
is easy to see that
\begin{align}
\begin{split}
& \{Ch(\mathcal{E})Td(T_{\mathcal{X}\backslash\mathcal{H}})\}_{(n+1,n+1)}
= \left((n+1)c_{1}(\mathcal{K}^{-1})c_{1}({\mathcal L})^{n}-\mu c_{1}({\mathcal L})^{n+1}\right)\\ \\
\end{split}
\end{align}
All we need to do is calculate the first Chern class of
$\textbf{L}(m)$. This is the content of the following
\begin{proposition}
\begin{align}
\begin{split}
\{Ch(\mathcal{E}({m})Td(T_{\mathcal{X}\backslash\mathcal{H}})\}_{(n+1,n+1)}=
(n+1)c_{1}(\mathcal{K}^{-1})c_{1}({\mathcal L})^{n}-\mu
c_{1}({\mathcal L})^{n+1}
\end{split}
\end{align}
\end{proposition}
The proof is straightforward. To begin with, we have
\begin{align}
\begin{split}
Ch(\mathcal{E}({m}))=&
(-1)^{n}2(n+1)\{\sum_{i=0}^{n}(-1)^{i}\binom{n}{i}\sum_{j=0}^{\infty}\frac{(m+i)^{j}}{j!}c_{1}({\mathcal L})^{j}\\
&+m
\sum_{i=0}^{n+1}(-1)^{i}\binom{n+1}{i}\sum_{j=0}^{\infty}\frac{(m+i)^{j}}{j!}c_{1}({\mathcal
L})^{j}
\}\\
&-(-1)^{n+1}(\mu
+n(n+1))\sum_{i=0}^{n+1}(-1)^{i}\binom{n+1}{i}\sum_{j=0}^{\infty}\frac{(m+i)^{j}}{j!}
c_{1}({\mathcal L})^{j}
\end{split}
\end{align}
Now switch the order of summation and use the binomial identities (ignore powers higher than $n+2$) to get
\begin{align*}
\begin{split}
\sum_{j=0}^{\infty}\sum_{i=0}^{n}(-1)^{i}\binom{n}{i}\frac{(m+i)^{j}}{j!}c_{1}({\mathcal
L})^{j} =&
\sum_{i=0}^{n}(-1)^{i}\binom{n}{i}\frac{(m+i)^{n}}{n!}c_{1}({\mathcal L})^{n}\\
&+\sum_{i=0}^{n}(-1)^{i}\binom{n}{i}\frac{(m+i)^{n+1}}{(n+1)!}c_{1}({\mathcal L})^{n+1}\\
=& (-1)^{n}c_{1}({\mathcal L})^{n}+  (-1)^{n}(m+\frac{n}{2})c_{1}({\mathcal L})^{n+1}\\
=& (-1)^{n}c_{1}({\mathcal L})^{n}+ (-1)^{n}mc_{1}({\mathcal
L})^{n+1}+ (-1)^{n}\frac{n}{2}c_{1}({\mathcal L})^{n+1}
\end{split}
\end{align*}
Therefore
\begin{align*}
\begin{split}
&Ch(\mathcal{E}(m))=\\
&(-1)^{n}2(n+1)\{(-1)^{n}c_{1}({\mathcal L})^{n}+
(-1)^{n}\frac{n}{2}c_{1}({\mathcal L})^{n+1}\} -(\mu
+n(n+1))c_{1}({\mathcal L})^{n+1}
\end{split}
\end{align*}
So that
\begin{align*}
\begin{split}
&\{Ch(\mathcal{E}(m))Td(T_{\mathfrak{X}\backslash \mathcal{H}})\}_{(n+1,n+1)}=\\
&\left(2(n+1)c_{1}({\mathcal L})^{n}c_{1} - \mu c_{1}({\mathcal
L})^{n+1}\right)(1-\frac{1}{2}c_{1}
({\mathcal K}_{\mathfrak{X}\backslash \mathcal{H}}))\\
&= -(n+1)c_{1}({\mathcal L})^{n}c_{1}({\mathcal
K}_{\mathfrak{X}\backslash \mathcal{H}})-\mu c_{1}({\mathcal
L})^{n+1}.
\end{split}
\end{align*}
So the proposition is proved. This proposition tells us that the
determinant lines have the {same} first Chern class. Therefore, they are (equivariantly)
isomorphic.

\section{The Weight of the CM polarization}

Since our two linearisations
are isomorphic they have the \emph{same weight} under the action
of any $\lambda:\mathbb{C}^{*}\rightarrow \gc$. That is
\begin{align}
w_{\lambda}(\textbf{L}(m))=
w_{\lambda}(\textbf{L}_{CM}^{-1}) 
\end{align}
 Again we stress that the weight on the left hand side  
is \emph{independent} of
the parameter $m$ that appears in its definition, and we are invited 
to let this parameter tend to infinity. 
This is analogous to the heat equation proof of the index theorem of Atiyah 
and Singer. In that case one has an alternating \emph{sum} of \emph{traces} where the parameter tends to \emph{zero}, in our case we have an alternating \emph{product} of \emph{determinants} and the parameter tends to \emph{infinity}.\\

This observation allows us to express the weight
$w_{\lambda}(\textbf{L}_{CM}^{-1})$ in terms of the
leading and subdominant coefficients
in the expansion of the weight of the mth Hilbert point
$w_{\lambda}(\mbox{Hilb}_{m}(X))$ of $X\subset \cpn$ under the action of a
1-psg. $\lambda$ of $\gc$.

First let us recall the definition of Hilbert points. The basic idea is simple enough: a projective variety is the zero set of finitely many homogenous polynomials, the issue is to consider polynomials \emph{which have the same degree}.

Let $X$ be a projective variety in $ \cpn$. Then we have for
$m>>0$ the exact sequence

\[ \begin{CD}
0@>>>H^{0}(\mathcal{I}_{X}(m))@>i>>H^{0}(\cpn,  {\mathcal
O}(m))@>Res_{m}(X)>> H^{0}(X,  {\mathcal O}(m)_{X}) @>>>0
\end{CD}
\]
where the vector space on the left is
\begin{align*}
\begin{split}
H^{0}(\mathcal{I}_{X}(m)))= &\{\mbox{All homogeneous polynomials of degree }\ m\\
&\mbox{in}\ N+1 \ \mbox{variables that vanish on}\ X\}
\end{split}
\end{align*}
If we let $P(m)$ be the Hilbert polynomial of $X$, then we have
\begin{align*}
\begin{split}
&H^{0}(\mathcal{I}_{X}(m))\in Gr(P(m),H^{0}(\cpn,  {\mathcal O}(m))) \\ \\
\end{split}
\end{align*}
where $Gr(P(m),H^{0}(\cpn,  {\mathcal O}(m)))$ denotes the
Grassmannian of \emph{codimension} $P(m)$ subspaces of
$H^{0}(\cpn,  {\mathcal O}(m))$. Using the Plucker embedding, we
may associate to $X\subset \cpn$ the point in the following
projective space
\begin{align*}
 {\bf{det}}(H^{0}(\mathcal{I}_{X}(m)))\in
 \mathbb{P}(\bigwedge^{\binom{N+m}{m}-P(m)}H^{0}(\cpn,  {\mathcal O}(m))).
\end{align*}
The $\mbox{m}^{th}$ Hilbert point of $X$ with respect to the given
polarization ${\mathcal O}(1)|_X$ is given by its dual
${\bf{det}}(H^{0}(X,  {\mathcal O}(m)_{X})$. To fix notation, we
will denote this Hilbert point by
\begin{align*}
\mbox{Hilb}_{m}(X):= {\bf{det}}(H^{0}(X,  {\mathcal O}(m)_{X}) \in
\mathbb{P}(\bigwedge^{P(m)}H^{0}(\cpn,  {\mathcal O}(m))^{*}).
\end{align*}
Since $\gc$ acts on this big projective space, we can associate a
weight to each 1psg. $\lambda:\mathbb{C}^{*}\rightarrow \gc$. Namely the weight of the action on $\mbox{Hilb}_{m}(X)$.  It is easy
to see that any such action can be diagonalized on $E$. Let
$\{e_{1},\dots,e_{d}\}$ be such a basis, i.e.
$\lambda(\alpha)e_{i}= \alpha^{m_{i}}e_{i} \quad \left(m_{i}\in
\mathbb{Z}\right)$. Next express any $v\in E$ in terms of this
basis $v =\sum_{i=1}^{d}v_{i}e_{i} $. Then the \emph{slope} of
$\mbox{Hilb}_{m}(X)$ relative to $\lambda$ is the number (usually
denoted by $\mu(\lambda,v)$:
\begin{align}
\mbox{Max}\{-m_{i}|v_{i}\neq 0\}
\end{align}

We define the weight $w_{\lambda}( \mbox{Hilb}_{m}(X))$ of
$\mbox{Hilb}_{m}(X)$ to be -$\mu(\lambda,v)$ for any $v$
lifting ${\bf{det}}(H^{0}(X,  {\mathcal O}(m)_{X}))$.

Our aim is to study the weight $w_{\lambda}(\mbox{Hilb}_{m}(X))$.
This is given by a numerical
polynomial of degree at most $n+1$ where $n=\mbox{dim}(X)$. In
other words
\begin{align*}
w_{\lambda}(\mbox{Hilb}_{m}(X))= a_{n+1}m^{n+1}+
a_{n}m^{n}+O(m^{n-1}).
\end{align*}
Recall that $\textbf{L}(m)$ was defined to be the determinant of
the total direct image of the virtual bundle\footnote{${\mathcal
L}$ is the pullback of ${\mathcal O}_{\cpn}(1)$}
\begin{align}
\begin{split}
\mathcal{E}(m):= & (-1)^{n}2(n+1)\{{\mathcal L}^{m}({\mathcal O}-
{\mathcal L})^{n}+m{\mathcal L}^{m}({\mathcal O}-{\mathcal L})^{n+1}\}\\
&-(-1)^{n+1}(\mu+n(n+1)){\mathcal L}^{m}({\mathcal O}-{\mathcal
L})^{n+1}
\end{split}
\end{align}

Observe that, since $m>>0$ all the higher direct image sheaves vanish
\begin{align*}
{R^{i}_{\pi}}_{*}( \mathcal{L}^{(m+j)})=0 \qquad i>0,
\end{align*}
we have the canonical isomorphism of determinant lines for $z \in
\mathcal{H}$.
\begin{align*}
\textbf{L}(m)_{z}\equiv (\textbf{L}_{1}(m)\otimes \textbf{L}_{2}^{m}(m))^{(-1)^{n}2(n+1)}_{z}\otimes (\textbf{L}_{2}(m))^{
(-1)^{n}(\mu +n(n+1))}_{z}
\end{align*}
Where we have defined
\begin{align*}
\textbf{L}_{1}(m)_{z}=\bigotimes_{i=0}^{n}{\bf{det}}H^{0}(X_{z},\mathcal{O}(m+i))^{(-1)^{i}\binom{n}{i}}
\end{align*}

\begin{align*}
\textbf{L}_{2}(m)_{z}=\bigotimes_{i=0}^{n+1}{\bf{det}}H^{0}(X_{z},\mathcal{O}(m+i))^{(-1)^{i}\binom{n+1}{i}}
\end{align*}
Therefore, if $z=z^{\lambda(0)}$ the weight of the action of $\lambda$ on the line is given by
\begin{align*}
\begin{split}
w_{\lambda}(\mathcal{L}(m);z)=& (-1)^{n}2(n+1)w_{\lambda}(\textbf{L}_{1}(m),z)+(-1)^{n}2(n+1)mw_{\lambda}(\textbf{L}_{2}(m),z)\\&+ (-1)^{n}(\mu +n(n+1))w_{\lambda}(\textbf{L}_{2}(m),z)
\end{split}
\end{align*}
The weight of the actions on the $\textbf{L}_{i}(m)$ are given in terms of hilbert points. For example, the weight of the action on $\textbf{L}_{1}(m)$ is given by
\begin{align*}
\begin{split}
w_{\lambda}(\textbf{L}_{1}(m),z)=& \sum_{i=0}^{n}(-1)^{i}\binom{n}{i}w_{\lambda}(\mbox{Hilb}_{m+i}(X_{z}))\\
&\sum_{i=0}^{n}(-1)^{i}\binom{n}{i}(a_{n+1}(m+i)^{n+1}+a_{n}(m+i)^{n}+O((m+i)^{n-1}))\\
&=(-1)^{n}(m+\frac{n}{2})(n+1)!a_{n+1}+(-1)^{n}n!a_{n}
\end{split}
\end{align*}
Where we have used the binomial identities again. Of course, one does the same thing to calclulate the other weight $w_{\lambda}(\textbf{L}_{2}(m))$.
Now just put all these computations together to see that
\begin{align*}
w_{\lambda}(\textbf{L}_{CM}^{-1})=    w_{\lambda}(\textbf{L}(m))=(n+1)!\left(2a_{n}(\lambda)-\mu a_{n+1}(\lambda)\right)
\end{align*}
Let $P(m)$ be the Hilbert polynomial. Then for $m>>0$ we have
\begin{align*}
mP(m)= mh^{0}(X, {\mathcal O}(m))=
b_{n+1}m^{n+1}+b_{n}m^{n}+O(m^{n})
\end{align*}
Where the $b_{i}$ are given by Hirzebruch Riemann-Roch. Following
Donaldson, we let $F_{1}$ be the coefficient of $\frac{1}{m}$ in
the expansion below
\begin{align*}
\begin{split}
\frac{w_{\lambda}(Hilb_{m}(X))}{mP(m)}&= \frac{a_{n+1}m^{n+1}+ a_{n}m^{n}+O(m^{n-1})}{b_{n+1}m^{n+1}+m^{n}+O(m^{n})}\\
&= \frac{a_{n+1}}{b_{n+1}}+
\frac{a_{n}b_{n+1}-a_{n+1}b_{n}}{b_{n+1}^{2}}\frac{1}{m}+
O(\frac{1}{m^{2}})
\end{split}
\end{align*}
A simple computation shows that

\begin{align*}
F_{1}(\lambda, X):= \frac{n!}{2d}(2a_{n}-\mu a_{n+1})
\end{align*}
Theorem 1.2 follows.

In [Don02] Donaldson defines  the generalised Futaki
invariant of the degeneration $\lambda$ to be this $F_{1}$.
In that paper he observed that this coincides with the definition
of Tian \emph{when the central fiber is smooth}. Our results in
this paper refine this observation. Theorem (1.3) shows that the
weight $F_{1}$, is the leading term in the asymptotics of the
\emph{Reduced} K-energy map for an \emph{arbitrary} central fiber.
\section{The Reduced K Energy map and $F_{1}$}

In this section we relate the general algebraic results of the previous sections to the K-energy map. Everything we need has already appeared in [Tian94] and we refer the reader to that paper for more details.
In particular we show that  the leading term of the reduced K-Energy map is just the weight $F_{1}$ up to a positive multiple.

Let $\varphi_{t}$ be a smooth path in $ P(X,\omega)$ (all K\"ahler potentials) joining $0$ with $\varphi$.
Then the \emph{K-energy map} introduced in [BM] is given by
\begin{align}
 \qquad \nu_{\omega}(\varphi):= -\frac{1}{V}\int_{0}^{1}\int_{X}\dot{\varphi_{t}}(\mbox{Scal}(\varphi_{t})-\mu)\omega_{t}^{n}dt
\end{align}

$\mbox{Scal}(\varphi_{t})$ denotes the scalar curvature of the metric $\omega +\sqrt{-1}\dl\dlb\varphi_{t}$.


Now we recall how the K-energy map can be viewed as a \emph{norm} on the CM polarisation. More precisely, $\nu_{\omega}$ can be veiwed as the logarithm of a singular norm on the CM polarisation. This fact allows us to find the precise asymptotics of the (reduced) K-energy map along any 1psg.$\lambda$. For the moment we assume that $\mathfrak{X}$ is smooth.
Let $\eta$ be a smooth test form on $\gc$ of type (g-1, g-1) where g is the dimension of $\gc$.
Define $G^{\mathbb{C}}X_{z}:= \{(\sigma,y)\in \gc \times \cpn:y\in \sigma X_{z}\}$, we note that
$G^{\mathbb{C}}X_{z}$ is biholomorphic to $G^{\mathbb{C}}\times X_{z}$.
Then we have the following
\begin{proposition}\label{ken}[Tian94] (complex hessian of the K-energy map)\\
Let $\pi_{1}^{-1}(z)=X_{z}\subset \cpn$, where $z\in \mathcal{H}_{\infty}$.
\begin{align}
d\int_{G^{\mathbb{C}}}\nu_{\omega,z}(\varphi_{\sigma} )\dl\dlb \eta = \int_{G^{\mathbb{C}}X_{z}}(R_{G^{\mathbb{C}}|X_{z}}+\frac{\mu(X)}{n+1}p_{2}^{*}(\omega_{FS}))\wedge p_{2}^{*}(\omega_{FS}^{n})\wedge p_{1}^{*}\eta
\end{align}
\end{proposition}
Above $\mathcal{H}_{\infty}$ denotes the locus of  smooth fibers in $\mathcal{H}$, i.e. $X_{z}:= \pi_{1}^{-1}(z)$ is a
smooth subvariety of $\cpn$ for $z \in \mathcal{H}_{\infty}$.
$\nu_{\omega,z}$ denotes the K-energy map on $X_{z}$.

Let us give an explanation of the curvature term $R_{G^{\mathbb{C}}|X_{z}}$.

Observe that $\pi_{2}^{*}\omega_{FS}$ induces a K\"ahler  metric on $\pi_{1}^{-1}(z)$ $(z \in \mathcal{H}_{\infty})$ and hence a metric on the \emph{relative} canonical bundle $K_{X_{z}}$ which we denote by $R(\pi_{2}^{*}(\omega_{FS}))$.
Now let $g_{\mathfrak{X}}$ and $g_{\mathcal{H}}$ denote two K\"ahler metrics on
$\mathfrak{X}$ and $\mathcal{H}$ respectively.
In this way we obtain \emph{another} metric on the relative canonical bundle
\begin{align*}
\mathcal{K}_{\mathfrak{X}}:= K_{\mathfrak{X}}\otimes \pi_{1}^{*} K_{\mathcal{H}}^{-1}\equiv K_{X_{z}}
\end{align*}
over the smooth locus. We let $R_{\mathfrak{X}|\mathcal{H}}$ denote its curvature.
\begin{align*}
R_{\mathfrak{X}|\mathcal{H}}:= R(g_{\mathfrak{X}})-\pi_{1}^{*}R(g_{\mathcal{H}})
\end{align*}

In the diagram below $p_{z}$ denotes the evaluation map, i.e.
$p_{z}(\sigma):= \sigma z$
\[ \begin{CD}
p_{z}^{*}\mathfrak{X}=\gc X_{z}@>p_{z,2}>> \mathfrak{X}@>\pi_{2}>>\cpn \qquad (**)\\
@Vp_{z,1}VV @V\pi_{1}VV \\
\gc@>p_{z}>>\mathcal{H}\\
\end{CD}
\]
Then
\begin{align*}
R_{\gc |X_{z}}:= p_{z,2}^{*}(R(\pi_{2}^{*}(\omega_{FS})))
\end{align*}

The relationship between these two choices is given by the next proposition.
\begin{proposition}\label{comp}
There is a smooth function $\Psi$ defined away from $\mathfrak{X}_{sing.}$
\begin{align*}
\Psi:\mathfrak{X}\backslash \pi_{1}^{-1}(\Delta)\rightarrow \mathbb{R}
\end{align*}
such that
\begin{align*}
 R(\pi_{2}^{*}(\omega_{FS}))=R_{\mathfrak{X}}-\pi_{1}^{*}(R_{\mathcal{H}})+\sqrt{-1}\dl\dlb \Psi
\end{align*}
Moreover, if we define $\Psi_{\mathcal{H}}(z):= \int_{\{y\in\pi_{1}^{-1}(z)\}}\Psi(y) \pi_{2}^{*}(\omega_{FS})^{n}$. Then $\Psi_{\mathcal{H}}(z)$ is bounded from above, is smooth outside $\Delta$ (the locus of singular fibers), continuous outside $\Delta_{{m}}$, and goes to $-\infty$ as $z\rightarrow \Delta_{{m}}$. $\Delta_{{m}}$ denotes the locus of $z\in \mathcal{H}$ where $X_{z}$ has a component of multiplicity greater than one.
\end{proposition}

\begin{proposition}\label{c0}
There is a continuous hermitian metric $|| \ ||_{CM}$ on $\textbf{L}_{CM}^{-1}
$ such that
\begin{align*}
-R(|| \ ||_{CM})= {\pi_{1}}_{*}\left((n+1)R_{\mathfrak{X}|\mathcal{H}}\wedge \pi_{2}^{*}(\omega_{FS})^{n}+\mu\pi_{2}^{*}(\omega_{FS})^{n+1}\right)
\end{align*}
In the weak sense.
\end{proposition} 

Recall that for $\sigma \in \gc$ we define $\varphi_{\sigma}$ by the relation
\begin{align*}
\sigma^{*}\omega_{FS}= \omega_{FS}+\dl\dlb \varphi_{\sigma}
\end{align*}

Below
$\nu_{\omega,z}(\sigma)$ denotes the K energy of $X_{z}$ applied to the potential $\varphi_{\sigma}$.
\begin{proposition}
For every smooth test form $\eta$ on $\gc$ of type (g-1,g-1), we have
\begin{align*}
d(n+1)\int_{\gc}\nu_{\omega,z}(\sigma)\dl\dlb \eta = \int_{\gc}{\log}\left(e^{(n+1)\Psi_{\mathcal{H}}(\sigma z)}\frac{||\ ||_{CM}^{2}(\sigma z)}{||\ ||_{CM}^{2}( z)}\right)\dl\dlb \eta
\end{align*}
\end{proposition}
Putting everything together, and using the fact that $\pi_{1}(\gc)=1$, we have
\begin{theorem}\label{sing}
Let $z\in \mathcal{H}_{\infty}$, then
\begin{align*}
d(n+1)\nu_{\omega,z}(\sigma)= {\log}\left(e^{(n+1)\Psi_{\mathcal{H}}(\sigma z)}\frac{||\ ||_{CM}^{2}(\sigma z)}{||\ ||_{CM}^{2}( z)}\right)
\end{align*}
\end{theorem}
Now, let $\lambda$ be an algebraic one parameter subgroup of $\gc$, and let $z\in \mathcal{H}$. Let $z^{\lambda(0)}$ denote the limit of $z$ under this action.
Then ${\textbf{L}_{CM}^{-1}}|{z^{\lambda(0)}}$ is a one dimensional representation of $\lambda:\mathbb{C}^{*}\rightarrow \gc$. This $\mathbb{C}^{*}$ acts via a character $w_{\lambda}(\textbf{L}_{CM}^{-1},\ z)\in \mathbb{Z}$
$$
\lambda(\alpha)v=\alpha^{w_{\lambda}(\textbf{L}_{CM}^{-1},\ z)}v \qquad v\in {\textbf{L}_{CM}}^{-1}|{z^{\lambda(0)}}.
$$
Applying Theorem (\ref{sing}) gives the asymptotics of $\nu_{\omega,z}(\lambda(t))$ as $t\rightarrow 0$

\begin{align}
d\nu_{\omega,z}(\lambda(t))- \Psi_{\mathcal{H}}(z^{\lambda(t)}) = 2w_{\lambda}(\textbf{L}_{CM}^{-1},\ z)\mbox{log}(t)+O(1)
\end{align}
Now we turn our attention to the case at hand.
Recall that $\mathfrak{X}:= \overline{G^{\mathbb{C}}X}$ the closure of 
$G^{\mathbb{C}}X$ inside of $\overline{G^{\mathbb{C}}}\times \cpn$.
Then $\mathfrak{X}$ has a divisor singularity. Let $\mathfrak{X}_{\infty}$ be a resolution of singularities of $\mathfrak{X}$. Let $\Delta_{i}\ 1\leq i \leq k$ be the exceptional divisors. We remark that $\mathfrak{X}_{\infty}\backslash \sum_{ 1\leq i \leq k}\Delta_{i}$ is just $G^{\mathbb{C}}X$. Let $p_{2}$ be the obvious map from 
$\mathfrak{X}_{\infty}$ to $\cpn$ (essentially the second projection), then GRR yields at once that
\begin{align*}
c_{1}(\textbf{L}_{CM})= -c_{1}(\textbf{Det}_{\mathcal{H}}(R^{\bull}_{\pi_{\infty}*}p_{2}^{*}\mathcal{E}_{CM}))+(n+1) \sum_{ 1\leq i \leq k}c_{1}(\textbf{L}_{i})
\end{align*}
Here $\pi_{\infty}$ is the map from $\mathfrak{X}_{\infty}$ onto $\mathcal{H}$,
and 
\begin{align*}
\textbf{L}_{i}:= \textbf{Det}_{\mathcal{H}}(R^{\bull}_{\pi_{\infty}*}(\mathcal{O}(\Delta_{i})-\mathcal{O})(p_{2}^{*}\mathcal{O}_{\mathbb{P}^{N}}(1) -\mathcal{O} )^{n})
\end{align*}
Let
\begin{align*}
i_{X}:\gc\rightarrow \mathcal{H}
\end{align*}
denote the inclusion map.
Now equip everything with metrics and apply essentially the same argument leading to the proof of Theorem 4.1. What results is the following proposition
\begin{proposition}
Let $\eta$ be a smooth compactly supported (g-1,g-1) form on $\gc$.
Then
\begin{align*}
-\int_{\gc}i_{X}^{*}c_{1}(\textbf{L}_{CM}^{-1})\wedge \eta=(n+1)\int_{\gc}\dl\dlb \left(d\nu_{\omega}(\sigma)
-\Psi_{\mathcal{H}}(\sigma)-\sum_{1\leq i \leq k}\theta_{i}(\sigma)\right)\wedge \eta
\end{align*}
\end{proposition}
In this formula $\Psi$ is the same function as in \ref{comp}. The $\theta_{i}$ are defined as follows.
Let $S_{\Delta_{i}}$ be a defining section of $\mathcal{O}(\Delta_{i})$. This line bundle has been equipped with some hermitian metric. We set $\theta_{i}$ to be the integral
\begin{align*}
\int_{\sigma X}\log||S_{\Delta_{i}}||^{2}c_{1}(L)^{n} 
\end{align*}
This leads at once to 
\begin{align*}
\frac{1}{n+1}{\log}\left(\frac{||\ ||_{CM}^{2}(\sigma)}{||\ ||_{CM}^{2}(\mbox{e})}\right)= d\nu_{\omega}(\sigma)-\Psi_{\mathcal{H}}(\sigma)-\sum_{1\leq i \leq k}\theta_{i}(\sigma)
\end{align*}

Theorem 1.3 follows if we just set $\Psi_{X}(\sigma):= \Psi_{\mathcal{H}}(\sigma)+\sum_{1\leq i \leq k}\theta_{i}(\sigma)$.\\ \\


We would like to end this note with the following questions.

\begin{tabular}{ll}
\\
\\
\bull If a manifold is K-Stable, does it have stable Chow (and Hilbert points)?
\\
\\
\bull Suppose a 1psg K-destabilises $X$ (so $F_{1}>0$), does there exist \emph{another} 1psg\\ that also
destabilises $X$ with \emph{normal} limit cycle $X^{\lambda(0)}$?
\\
\\

\bull Is the CM polarisation \emph{ample} on the compliment of a proper subvariety of $\mathcal{H}$?\\ (i.e does the K-stability have a G.I.T. interpretation?)
\end{tabular}

\providecommand{\bysame}{\leavevmode\hbox to3em{\hrulefill}\thinspace}

 \small{\textsc{Department of Mathematics, Columbia University  NY. NY. 10027 }}\\
 {\em{E-mail:}} 
 {\texttt{stpaul@math.columbia.edu}}
 \\
 \small{\textsc{Department of Mathematics, MIT,
 Cambridge, MA 02139}}\\
 {\em{E-mail:}} 
 {\texttt{tian@math.mit.edu}}

\end{document}